\documentclass{article}

\usepackage[margin=1in]{geometry}
\usepackage{amsmath}
\usepackage{amssymb}
\usepackage{enumitem}

\title{A Proof of the Conjecture by Carpentier-De Sole-Kac}
\author{Keaton Stubis}
\date{October 20, 2014}

\begin{document}
\maketitle

\begin{center}
0. Introduction
\end{center}

Let $\mathcal{K}$ be a differential field with derivation $\partial$, and let $R$ be a differential subring of $\mathcal{K}$. We consider an $n \times n$ matrix $M$, whose elements lie in the ring $R[\partial]$. A number of useful notions can be defined on such a matrix, such as the following ones:
\newline\newline
For a matrix $A$ with entries in the field $\mathcal{K}[\partial]$, the Dieudonn\'e Determinant has the form $det A = det_1 A\lambda^d$, where $det_1 A \in \mathcal{K}$, $\lambda$ is indeterminate, and $d$ is an integer. Some of its characterizing properties are $det A \cdot det B = det AB$, that $det A$ changes sign upon permuting two rows of $A$, and that subtracting $h$ times $A$'s ith row from its jth row leaves $det A$ unchanged, when $h \in \mathcal{K}[\partial]$ and $i \neq j$. Furthermore, if $A$ is upper triangular, then $det_1 A$ is the product of the leading coefficients of the diagonal entries, and $d$ is the sum of their orders. In such a case, if any of the diagonal entries are 0, then these evaluate to $det_1 A=0$ and $d=-\infty$. That a determinant with these properties exists is shown in a more general context in [Die43].
\newline\newline
\textbf{Definition 0.1.} The total order of a matrix $A$ is
$$
tord(A) = \max_{\sigma \in S_n}\sum_{i=1}^{n}ord(A_{i,\sigma (i)})
$$
where $S_n$ denotes the group of permutations of $\{1,2,\dots n\}$. We can then also define the degeneracy degree of $A$ by
$$
dd(A) = tord(A) - d(A)
$$
where $d(A)$ is the order of $det A$.
\newline\newline
\textbf{Definition 0.2.} A system of integers $(N_1,N_2,\dots N_n,h_1,h_2,\dots,h_n)$ is called a majorant of $A$ if for all $i,j \in \{1,2,\dots n\}$
$$
ord(A_{ij})\leq N_j - h_i
$$
Given a matrix A and a majorant of that matrix, one can associate a characteristic matrix $\bar{A}(\lambda)$ to the majorant by the equation
$$
\bar{A}_{ij}(\lambda) = A_{ij;N_j-h_i} \lambda^{N_j - h_i}
$$
where $A_{ij;N_j-h_i}$ is the coefficient of $\partial^{N_j-h_i}$ in $A_{ij}$. We can also note that fact that for any majorant
$$
\sum_{i=1}^{n}N_i-h_i \geq tord(A)
$$
We call a majorant optimal when equality holds in the above statement.

\newpage

The following theorem of [CDSK12] follows from the results in [Huf65].
\newline\newline
\textbf{Theorem 0.3} Let $A$ be a matrix with elements in $\mathcal{K}[\partial]$ and let $det A \neq 0$. Then:
\begin{enumerate}[label = \roman*.]
\item $dd(A) \geq 0$

\item there exists an optimal majorant of A

\item if $dd(A) \geq 1$, then $det(\bar{A}(\lambda)) = 0$ for any majorant

\item if $dd(A) = 0$, then $det(\bar{A}(\lambda)) = 0$ for any majorant which is not optimal, and $det(\bar{A}(\lambda)) = det A$ for any majorant which is optimal
\end{enumerate}

It is obvious that $det_1 A \in R$ if $dd(A) = 0$, and it was shown in [CDSK12] that $det_1 A$ always lies in the integral closure of $R$ in $K$. However, in general $det_1 A \notin R$; there is a counterexample in [CDSK12] with $dd(A) = 2$. It was conjectured in [CDSK12] that $det_1 A \in R$ when $dd(A) = 1$. In the present paper this conjecture is proved.

\begin{center}
1. Proof of the Carpentier-De Sole-Kac conjecture.
\end{center}

\textbf{Proposition 1.1.} Let $\mathcal{K}$ be a differential field, and let $R$ be a differential subring that is a subring of $\mathcal{K}$. Now, let $M$ be a matrix whose elements lie in $R[\partial]$, and let $D = det_1 M$. If $dd(M) = 1$, then $D \in R$.
\newline\newline
Proof. First, note that we may multiply we may increase the order of each term in a given row or column by either left or right multiplying by a matrix of the form $diag(1\: \ldots 1\: \partial \: 1\: \ldots 1)$.
This operation results in a matrix whose entries are still in $R$ and also leaves the coefficient of the determinant unchanged. Furthermore, it increases both the total order and degree of the determinant by 1, leaving the degeneracy degree unchanged. Next, we would like to extract an optimal majorant $N_1, N_2,\ldots N_n, h_1,h_2,\ldots h_n$. Both a definition and some of its basic properties can be found in [CDSK12,Def4.6] and [CDSK,Thm4.7]. A majorant exists as long as $D\neq 0$, which we may assume, since otherwise $D=0$, at which point $D$ must lie in $R$ anyway. Now, define $N$ to be the largest of all the $N_i$ and $h$ the minimum of the $h_i$. Note that when we use the process described earlier to increase the degrees of the elements of column i by 1, we can get an optimal majorant for the new matrix by simply taking the old majorant and increasing $N_i$ by 1. Similarly, if we increase the degrees in row i by 1, then the new majorant is just the old majorant, with $h_i$ decreased by 1.

We can now use this to create a new matrix $M'$ as follows. Increase the degrees in column i by 1 $N - N_i$ times, and increase the degrees in row i by 1 $h_i - h$ times. Then, $M'$ will have a determinant coefficient of $D$, degeneracy degree 1, and an optimal majorant of $N,N,\ldots N,h,h,\ldots h$. From this optimal majorant, we may extract two properties of $M'$. First, every entry has degree less than or equal to $N-h$. Also, the total order is $n(N-h)$, found by summing over the majorant. These two facts together imply that every row must contain at least one term of degree exactly $N-h$. In fact, they show something stronger, but this is all we will need.

We now know that the largest degree term has degree $N-h$. We can now express the matrix $M'$ in the form $A\partial^{N-h} + B\partial^{N-h-1} \ldots$, where $A$ is the leading coefficient matrix of terms of degree $N-h$, $B$ contains the terms of degree $N-h-1$, and so on.

Now, we know that since $M'$ is degenerate, its characteristic matrix has determinant 0 (see [CDSK12,Thm4.7]). In this case, the characteristic matrix is exactly $A$. Letting the rows of $A$ be called $A_1,A_2,\ldots A_n$, the fact that $A$ has determinant 0 implies a linear relation amongst the rows with coefficients in the fraction field of $R$. We may clear denominators to get a relation with coefficients strictly in the ring $R$.
$$c_1A_1 + c_2A_2 \ldots + c_nA_n = 0$$

Next, multiply the first row of $M'$ by $c_1$ via the use of an elementary matrix. From the total order and degeneracy of $M'$, we can see that the determinant of this product is $c_1D\lambda^{n(N-h)-1}$. We can further alter the matrix by next subtracting $c_i$ times the ith row from the first row, for each $i\neq1$. This operation leaves the determinant unchanged, and yields a new matrix we call $M''$. We express $M''$ as a sum similar to the one we used before, so that $M'' = A'\partial^{N-h} + B'\partial^{N-h-1} \ldots$ for new matrices $A'$ and $B'$. These matrices take on the following forms:
$$
A'=
\begin{bmatrix}
0\\
A_2\\
A_3\\
\vdots\\
A_n\\
\end{bmatrix}
\qquad\qquad
B'=
\begin{bmatrix}
c_1B_1 + c_2B_2 \ldots + c_nB_n\\
B_2\\
B_3\\
\vdots\\
B_n\\
\end{bmatrix}
$$

We now note that the total order of $M''$ is at least $n(N-h)-1$, the order of its determinant. On the other hand, $N,N,\ldots N,h+1,h,h,\ldots h$ is a majorant, which implies that not only is it an optimal majorant, but that $M''$ is non-degenerate, and thus that its determinant is the determinant of the characteristic matrix [CDSK12,Thm.4.7]. To finish the proof, we now show that the determinant of the characteristic matrix is a multiple of $c_1$. As this is equal to $c_1D$, we can then apply the cancellation law to show that $D$ lies in the ring $R$.

The characteristic matrix is
$$
M''_{char} = 
\begin{bmatrix}
c_1B_1 + c_2B_2 \ldots + c_nB_n\\
A_2\\
A_3\\
\vdots\\
A_n\\
\end{bmatrix}
$$

Since the determinant is linear in the first row, we can set the determinant of $M''_{char}$ equal to $det(M''_1) + det(M''_2) + \ldots + det(M''_n)$, where
$$
M''_i = 
\begin{bmatrix}
c_iB_i\\
A_2\\
A_3\\
\vdots\\
A_n\\
\end{bmatrix}
$$

$M''_1$ has a first row that is a multiple of $c_1$, so its determinant is already a multiply of $c_1$. We now look at the other $M''_i$:
Since we can multiply by elementary matrices, we can divide out the factor of $c_i$ from the top row and multiply the ith row by $c_i$, leaving the determinant unaffected. Next, for each $j\neq1$, add $c_j$ times row j to row i. The resulting row i can then be simplified using our earlier linear relation and becomes a multiple of $c_1$. Thus, the determinant of $M''_i$ is a multiple of $c_1$, as desired. Below, these last steps are written out explicitly:
$$
\begin{bmatrix}
c_iB_i\\
A_2\\
A_3\\
\vdots\\
A_n\\
\end{bmatrix}
\rightarrow
\begin{bmatrix}
B_i\\
A_2\\
A_3\\
\vdots\\
c_iA_i\\
\vdots\\
A_n\\
\end{bmatrix}
\rightarrow 
\begin{bmatrix}
B_i\\
A_2\\
A_3\\
\vdots\\
c_2A_2+c_3A_3+\ldots+c_nA_n\\
\vdots\\
A_n\\
\end{bmatrix}
\rightarrow
\begin{bmatrix}
B_i\\
A_2\\
A_3\\
\vdots\\
-c_1A_1\\
\vdots\\
A_n\\
\end{bmatrix}
$$
\newline\newline
\begin{center}
References
\end{center}
\begin{flushleft}
[CDSK12] S. Carpentier, A. De Sole, V.G. Kac, \textit{Some Algebraic Properties of Differential Operators}, arXiv:1201.1992v1

[Die43] J. Dieudonn\'e \textit{Les d\'eterminants sur un corps non commutatif}, Bull. Soc. Math. France \textbf{71}, (1943), 27-45.

[Huf65] G. Hufford, \textit{On the characteristic matrix of a matrix of differential operators}, J. Differential Equations \textbf{1}, (1965) 27–38.
\end{flushleft}
\end{document}